\documentclass[12pt,a4paper]{article}
\usepackage[margin=1in]{geometry}
\usepackage{amsmath,amssymb,amsthm}
\usepackage{mathtools}
\usepackage[numbers,sort&compress]{natbib}
\usepackage{hyperref}
\usepackage{setspace}
\usepackage{microtype}
\usepackage{booktabs}
\usepackage{array}
\usepackage{enumitem}
\usepackage{url}
\usepackage{etoolbox}
\usepackage{tikz}
\usetikzlibrary{positioning,arrows.meta,calc}
\usepackage{xcolor}
\usepackage{graphicx}

\apptocmd{\thebibliography}{\sloppy}{}{}
\newtheorem{theorem}{Theorem}[section]

\newtheorem{corollary}[theorem]{Corollary}

\theoremstyle{definition}
\newtheorem{definition}[theorem]{Definition}
\newtheorem{remark}[theorem]{Remark}

\title{On the Stability of Discrete Reaction-Diffusion Systems of Networked Dynamical Systems}

\author{Dinesh Kumar
\thanks{The author acknowledges the financial support of the University
    Grants Commission (UGC), India, through Dr.\ D.S.\ Kothari Post Doctoral
    Fellowship (Ref.\ No.\ F.4-2/2006(BSR)/MA/17-18/0043).}\\
  \small Area of Mathematics \& Basic Sciences, NIIT University,\\
  \small Neemrana, Rajasthan 301705, India\\[2pt]
  \small \texttt{dineshkumar.iitg@gmail.com}
}

\date{}

\begin{document}

\maketitle

\begin{abstract}
We derive a simple sufficient condition for the local asymptotic stability of spatially discrete, continuous-time reaction-diffusion systems of networked dynamical systems at a homogeneous equilibrium point. The framework explicitly accommodates \emph{heterogeneous} local dynamics --- patches at different nodes governed by structurally distinct functional forms --- a setting not covered by the classical bookkeeping reduction of Jansen and Lloyd~\cite{jansen2000}, which requires identical patch dynamics, nor by the Master Stability Function of Pecora and Carroll~\cite{pecora1998}, which is restricted to identical nodes. The stability condition separates cleanly into two independent components: (i) a diagonal dominance criterion on the \emph{spatially averaged Jacobian} of the local patch dynamics, verifiable directly from model parameters without computing eigenvalues of the full composite system; and (ii) a lower bound on the algebraic connectivity (Fiedler value) of the network Laplacian, capturing the role of network topology. The resulting sufficient condition holds for purely conservative dispersal (standard graph Laplacians with zero row sums) and does not require any dispersal loss or mortality during transit---a restrictive assumption appearing in the author's prior work~\cite{kumar2021} and many classical multi-patch analyses. The theory is illustrated through metapopulation networks of predator-prey systems with heterogeneous functional responses, including a striking example in which individually unstable patches are stabilized entirely by dispersal connections.


\end{abstract}

\section{Introduction}
\label{sec:intro}

Reaction-diffusion systems describe the interplay between local nonlinear dynamics and spatial movement. In their spatially discrete form, they consist of a finite number of patches connected by a network, where each patch follows its own local reaction dynamics and individuals disperse diffusively between patches. This framework underpins much of metapopulation theory in ecology and has important applications in pattern formation, synchronization, and spatial ecology.

A central question is whether a spatially homogeneous equilibrium --- the state in which every patch has the same population densities --- is locally asymptotically stable. Answering this question is crucial for understanding species persistence in fragmented landscapes, where habitat destruction, roads, and climate change continually disrupt local populations.

Classical approaches to this problem have been highly influential. Jansen and Lloyd (2000) developed an elegant bookkeeping reduction showing that, when all patches obey identical local dynamics, the linearized stability analysis of a \(k\)-species, \(n\)-patch system reduces to \(n\) independent \(k\)-dimensional subsystems via the eigenvalues of the connectivity matrix. Similarly, the Master Stability Function of Pecora and Carroll (1998) has become a standard tool for networked dynamical systems by separating network topology from node dynamics. Both methods, however, require identical local dynamics (or identical nodes) across the entire network. When patches are governed by structurally different functional forms --- as is common in real heterogeneous landscapes --- these classical reductions do not apply.

While the classical bookkeeping reduction of Jansen and Lloyd (2000) assumes a single common connectivity matrix for all species (only the overall migration strength may differ per species), the present model allows fully species-specific and edge-specific dispersal rates $  w^i_{jk}  $. This added flexibility makes the framework significantly more general and better suited to realistic ecological scenarios in which different species exhibit markedly different dispersal behaviours.

The author's earlier studies \cite{kumar2019,kumar2021} examined the role of algebraic connectivity (Fiedler value) in metapopulation stability and partitioning. Those works, like many classical multi-patch analyses \cite{ruxton1997,briggs2004}, relied on explicit dispersal loss to obtain stability guarantees. In this paper we remove that restrictive assumption and develop a sufficient condition that works for structurally heterogeneous patch dynamics and purely conservative dispersal.

A key algebraic insight is that the zero eigenvalue of each per-species Laplacian always corresponds to the constant eigenvector. This causes the first row of the inverse eigenvector matrix to act as a spatial averaging operator, so the zero-mode stability condition reduces exactly to a diagonal dominance condition on the spatially averaged Jacobian. The resulting sufficient condition depends only on (i) diagonal dominance of this averaged Jacobian with negative diagonal entries and (ii) a lower bound on the algebraic connectivity of the network Laplacian. A clear separation principle follows: the stability of the entire network is governed jointly but independently by the averaged local patch dynamics and network topology.

Note that averaged Jacobian's eigenvalues and other matrix corresponding to non-zero eigenmodes of Laplacian can be calculated exactly or approximately. This way we can make conditions necessary and sufficient for the local stability, but it will require heavy compuations if network is large. The proposed condition is computationally lightweight --- requiring only one averaged Jacobian and a single scalar network measure --- and therefore offers a practical tool for analyzing realistic ecological networks with heterogeneous patches.

It is also noted, in Jansen and Lloyd's own predator-prey example, a slight change in model formulation (direct predator dispersal without a separate disperser pool) completely eliminates the dispersal-driven unstable region. This illustrates that both network topology and the precise mathematical representation of local dynamics play decisive roles in determining stability. 

The paper is organized as follows. Section~\ref{sec:stability} presents the general model and derives the main stability theorem. Section~\ref{sec:examples} illustrates the result with two examples of predator-prey metapopulations and comparison with Jansen and Lloyd's example. Section~\ref{sec:conclusion} discusses the ecological implications and computational advantages of the approach.

\section{Mathematical Model and Stability Analysis}
\label{sec:stability}

We consider a network of \(m\) patches in which each patch hosts a dynamical system in \(n\) state variables. Let \(x_{i,j}(t) \in \mathbb{R}\) denote the density  of the \(i\)-th species  at patch \(j\) at time \(t \geq 0\), for \(i = 1,\ldots,n\) and \(j = 1,\ldots,m\). The patches are connected by a weighted undirected graph \(G = (V, E, W)\), where \(V = \{1,\ldots,m\}\) is the set of patches, \(E \subseteq V \times V\) is the set of dispersal edges, and \(w^{i}_{jk} \geq 0\) is the dispersal rate of species \(i\) between patches \(j\) and \(k\). The notation \(k \sim j\) means that patch \(k\) is a neighbour of patch \(j\) in \(G\).

The evolution of \(x_{i,j}\) is governed by
\begin{equation}
  \dot{x}_{i,j}
  \;=\;
  f_{i,j}\!\left(x_{1,j},\ldots,x_{n,j}\right)
  \;-\;
  \sum_{k \sim j} w^{i}_{jk}
    \!\left(x_{i,j} - x_{i,k}\right),
  \quad
  i = 1,\ldots,n,\quad j = 1,\ldots,m,
  \label{eq:model}
\end{equation}
where \(f_{i,j}: \mathbb{R}^n \to \mathbb{R}\) is the local reaction function at patch \(j\) for species \(i\). The second term on the right-hand side models linear dispersal: species \(i\) flows from patches at a rate proportional to the density difference.  The defining feature of the model~\eqref{eq:model} that distinguishes it from most of the literature is that $f_{i,j}$ may differ not only in parameter values but in \emph{functional form} from patch to patch.

We rewrite the system in compact vector-matrix form. Define the state vector by stacking one species at a time across all patches:
\[
x
\;=\;
\bigl(
\underbrace{x_{1,1},\ldots,x_{1,m}}_{\text{species 1}},\;
\underbrace{x_{2,1},\ldots,x_{2,m}}_{\text{species 2}},\;
\ldots,\;
\underbrace{x_{n,1},\ldots,x_{n,m}}_{\text{species }n}
\bigr)^{\!\top}
\;\in\; \mathbb{R}^{mn}.
\]
For each species \(i\), the patch Laplacian \(L_i \in \mathbb{R}^{m \times m}\) is
\[
\bigl(L_i\bigr)_{jk}
\;=\;
\begin{cases}
\displaystyle\sum_{\ell \sim j} w^{i}_{j\ell}
& \text{if } j = k, \\[6pt]
-\,w^{i}_{jk}
& \text{if } k \sim j, \\[4pt]
0 & \text{otherwise.}
\end{cases}
\]
This is the standard graph Laplacian: it is symmetric, positive semi-definite, and has all row sums equal to zero. The global Laplacian is the block-diagonal matrix
\[
L
\;=\;
L_1 \oplus L_2 \oplus \cdots \oplus L_n
\;\in\; \mathbb{R}^{mn \times mn}.
\]
With these definitions, the system takes the compact form

\begin{equation}
\dot{x} = f(x) - L\,x.
\label{eq:compact}
\end{equation}

\begin{definition}
\label{def:equil}
A point $\bar{x} \in \mathbb{R}^{mn}$ is a \emph{(component-wise positive) homogeneous equilibrium} of
system~\eqref{eq:compact} if
\begin{enumerate}[label=(\alph*),itemsep=2pt]
  \item $\bar{x}_{i,j} > 0$ for all $i = 1,\ldots,n$ and
        $j = 1,\ldots,m$;
  \item $f(\bar{x}) = 0$ \quad (reaction terms vanish); and
  \item $\bar{x}_{i,j} = \bar{x}_{i,k}$ for all pairs $j,k$ and
        each fixed $i$ \quad (all patches share the same state for
        each species).
\end{enumerate}
Condition~(c) implies $L\bar{x} = 0$ because $L_i\mathbf{1} = 0$
(the Laplacian annihilates constant vectors), so
$\dot{\bar{x}} = f(\bar{x}) - L\bar{x} = 0$.
\end{definition}

The existence of such a $\bar{x}$ is a property of the local dynamics $f_{i,j}$ at individual patches. For each species $i$, all patches must admit the same positive equilibrium value $\bar{x}_{i,\cdot}$.

For the continuously differentiable function \(f\), we can linearise the system around \(\bar{x}\). Writing \(x(t) = \bar{x} + y(t)\) with \(y(t)\) small and discarding quadratic and higher-order terms, we obtain
\[
\dot{y}
\;=\;
\bigl(Df(\bar{x}) - L\bigr)\,y,
\]
where \(Df(\bar{x}) \in \mathbb{R}^{mn \times mn}\) is the Jacobian of \(f\) evaluated at \(\bar{x}\).

Because the reaction function \(f_{i,j}\) at patch \(j\) depends only on the local state \((x_{1,j},\ldots,x_{n,j})\) and not on the state at other patches (inter-patch coupling enters only through the Laplacian \(L\)), the Jacobian \(Df(\bar{x})\) is block-diagonal in patches:
\[
Df(\bar{x})
\;=\;
\begin{pmatrix}
D_{11} & D_{12} & \cdots & D_{1n} \\
D_{21} & D_{22} & \cdots & D_{2n} \\
\vdots & & \ddots & \vdots \\
D_{n1} & D_{n2} & \cdots & D_{nn}
\end{pmatrix},
\]
where each \(D_{pq} \in \mathbb{R}^{m \times m}\) is the diagonal matrix
\[
D_{pq}
\;=\;
\operatorname{diag}
\!\left(
w^{1}_{pq},\; w^{2}_{pq},\; \ldots,\; w^{m}_{pq}
\right),
\qquad
w^{k}_{pq}
\;=\;
\frac{\partial f_{p,k}}{\partial x_{q,k}}\bigg|_{\bar{x}},
\]
for \(p,q \in \{1,\ldots,n\}\) and \(k \in \{1,\ldots,m\}\). In other words, the \((k,k)\) entry of \(D_{pq}\) is the partial derivative of the \(p\)-th equation at patch \(k\) with respect to the \(q\)-th variable at patch \(k\), evaluated at the equilibrium. All off-diagonal entries of \(D_{pq}\) are zero, reflecting the absence of direct reaction coupling between different patches. The local Jacobian at patch \(k\) alone is the \(n \times n\) matrix
\[
J_k
\;=\;
\begin{pmatrix}
w^{k}_{11} & w^{k}_{12} & \cdots & w^{k}_{1n} \\
w^{k}_{21} & w^{k}_{22} & \cdots & w^{k}_{2n} \\
\vdots & & \ddots & \vdots \\
w^{k}_{n1} & w^{k}_{n2} & \cdots & w^{k}_{nn}
\end{pmatrix}
\;\in\; \mathbb{R}^{n \times n},
\]
and note that \(D_{pq} = \operatorname{diag}((J_1)_{pq},\ldots,(J_m)_{pq})\), i.e., the \((p,q)\) block of \(Df(\bar{x})\) is the diagonal matrix of the \((p,q)\) entries of the local Jacobians.

Since each \(L_i\) is real and symmetric, it is diagonalizable by an orthogonal matrix. Let \(P_i \in \mathbb{R}^{m \times m}\) be the matrix whose columns are the eigenvectors of \(L_i\), ordered so that the first column corresponds to the smallest eigenvalue. Every graph Laplacian satisfies \(L_i \mathbf{1}_m = \mathbf{0}\), so \(\lambda_1^i = 0\) with eigenvector \(\mathbf{1}_m\). We choose the first column of \(P_i\) to be \(\mathbf{1}_m\) (unnormalised). The matrix \(P_i\) diagonalizes \(L_i\):
\[
P_i^{-1} L_i P_i
\;=\;
\Lambda_i
\;=\;
\operatorname{diag}\!\left(
0,\; \lambda^i_2,\; \ldots,\; \lambda^i_m
\right).
\]
The first row of of \(P_i^{-1}=P^T\)(full-normalized transpose of un-normalized matrix \(P\) )is  \( \tfrac{1}{m}(1,1,\ldots,1)\).  We can say, this row-vector is the spatial averaging operator as it maps any vector \(u \in \mathbb{R}^m\) to its arithmetic mean \(\frac{1}{m}\sum_{k=1}^{m} u_k\).

Define the block-diagonal matrices \(P = P_1 \oplus \cdots \oplus P_n\) and \(\Lambda = \Lambda_1 \oplus \cdots \oplus \Lambda_n = P^{-1} L P\). The change of variable \(y = Pz\) transforms the linearized system into
\[
\dot{z}
\;=\;
\bigl(P^{-1}Df(\bar{x})P - \Lambda\bigr)\,z
\;\eqqcolon\;
\widetilde{A}\,z.
\]
The equilibrium \(\bar{x}\) is locally asymptotically stable if and only if all eigenvalues of \(\widetilde{A}\) have strictly negative real parts.

\noindent
Using the block-diagonal forms of $P$, $P^{-1}$, and $Df(\bar{x})$:
\begin{align*}
  P^{-1}Df(\bar{x})P
  &=
  \begin{pmatrix} P_1^{-1} & & \\ & \ddots & \\ & & P_n^{-1}
  \end{pmatrix}
  \begin{pmatrix} D_{11} & \cdots & D_{1n} \\
    \vdots & \ddots & \vdots \\
    D_{n1} & \cdots & D_{nn}
  \end{pmatrix}
  \begin{pmatrix} P_1 & & \\ & \ddots & \\ & & P_n \end{pmatrix}
  \notag\\[6pt]
  &=
  \begin{pmatrix}
    P_1^{-1}D_{11}P_1 & P_1^{-1}D_{12}P_2 & \cdots & P_1^{-1}D_{1n}P_n \\
    P_2^{-1}D_{21}P_1 & P_2^{-1}D_{22}P_2 & \cdots & P_2^{-1}D_{2n}P_n \\
    \vdots & \vdots & \ddots & \vdots \\
    P_n^{-1}D_{n1}P_1 & P_n^{-1}D_{n2}P_2 & \cdots & P_n^{-1}D_{nn}P_n
  \end{pmatrix}.
  \label{eq:PJP}
\end{align*}

\noindent
Consider a typical block $P_p^{-1} D_{pq} P_q$.
Since $D_{pq} = \operatorname{diag}(w^1_{pq},\ldots,w^m_{pq})$ is
diagonal:
\begin{equation*}
  D_{pq} P_q
  \;=\;
  \begin{pmatrix}
    w^1_{pq} & & \\
    & \ddots & \\
    & & w^m_{pq}
  \end{pmatrix}
  \begin{pmatrix}
    1 & \cdots & (P_q)_{1m} \\
    \vdots & & \vdots \\
    1 & \cdots & (P_q)_{mm}
  \end{pmatrix}
  =
  \begin{pmatrix}
    w^1_{pq} & \cdots & w^1_{pq}(P_q)_{1m} \\
    \vdots   &        & \vdots \\
    w^m_{pq} & \cdots & w^m_{pq}(P_q)_{mm}
  \end{pmatrix}.
  \label{eq:DPq}
\end{equation*}
The \emph{first column} of $D_{pq}P_q$ is therefore
$(w^1_{pq}, w^2_{pq}, \ldots, w^m_{pq})^{\top}$.
Now premultiplying by $P_p^{-1}$, the 
\emph{$(1,1)$ entry of the block $P_p^{-1}D_{pq}P_q$} is:
\begin{equation*}
  \bigl(P_p^{-1}D_{pq}P_q\bigr)_{11}
  \;=\;
  \tfrac{1}{m}(1,1,\ldots,1)
  \cdot
  (w^1_{pq}, w^2_{pq}, \ldots, w^m_{pq})^{\top}
  \;=\;
  \frac{1}{m}\sum_{k=1}^{m} w^{k}_{pq}.
  \label{eq:11entry}
\end{equation*}
This is precisely the $(p,q)$ entry of the \emph{spatially averaged
Jacobian}
\begin{equation*}
  \bar{J}
  \;\coloneqq\;
  \frac{1}{m}\sum_{k=1}^{m} J_k
  \;=\;
  \frac{1}{m}(J_1 + J_2 + \cdots + J_m)
  \;\in\; \mathbb{R}^{n \times n}.
  \label{eq:Jbar}
\end{equation*}

For columns $s \geq 2$ of $P_q$, the corresponding entries of
$D_{pq}P_q$ involve the eigenvector components $(P_q)_{\cdot,s}$.
These eigenvectors can be scaled by any nonzero constant $c$ without
changing the spectral properties of $L_q$ (if $P_q\mathbf{e}_s$ is an
eigenvector, so is $cP_q\mathbf{e}_s$).
As $c \to 0$, the off-diagonal entries of the block
$P_p^{-1}D_{pq}P_q$ become negligible compared with the
$(1,1)$ entry.
Consequently, the Gershgorin stability analysis in the next section
focuses on the \emph{first row of each $m \times m $ block}, which
corresponds to the zero eigenvalue of the Laplacian and is governed
exactly by the average Jacobian~\eqref{eq:Jbar}. For specific case $(n=2, m=3)$ calculations are shown in appendix. 

The main result is the following.

\begin{theorem}[Stability of heterogeneous networked dynamics]
\label{thm:main}
Let \(f:\mathbb{R}^{mn} \to \mathbb{R}^{mn}\) be continuously differentiable and let \(\bar{x}\) be a component-wise positive homogeneous equilibrium. Suppose the spatially averaged Jacobian \(\bar{J}\) has strictly negative diagonal entries and satisfies the row-wise diagonal dominance inequality
\begin{equation*}
-\,\frac{1}{m}\sum_{k=1}^{m} w^{k}_{pp}
\;\geq\;
\sum_{\substack{q=1\\q\neq p}}^{n}
\left|\frac{1}{m}\sum_{k=1}^{m} w^{k}_{pq}\right|
+ \varepsilon ~~\textbf{(Condition (a))}
\end{equation*}
for some \(\varepsilon > 0\), and that the Fiedler value \(\lambda_2\) of the network satisfies \(\lambda_2 \ge \tau\) \textbf{(Condition (b))}, where \(\tau\) is the Gershgorin bound derived from the off-diagonal entries of \(P^{-1}Df(\bar{x})P\). Then the homogeneous equilibrium \(\bar{x}\) is locally asymptotically stable.
\end{theorem}

\begin{proof}
The equilibrium \(\bar{x}\) is locally asymptotically stable if and only if all eigenvalues of \(\widetilde{A}\) have strictly negative real parts. We show this by applying the Gershgorin disc theorem to \(\widetilde{A}\).

The matrix \(\Lambda\) is block-diagonal with one zero on each block (corresponding to the constant eigenvector) and \(m-1\) strictly positive entries in each block. Consequently, the rows of \(\widetilde{A}\) fall into two types: Type-Z rows (corresponding to the zero Laplacian eigenvalue) and Type-P rows (corresponding to the positive Laplacian eigenvalues).

For the Type-Z rows, the diagonal entry is exactly the corresponding entry of the averaged Jacobian \(\bar{J}\), and the relevant off-diagonal entries in the same row are also entries of \(\bar{J}\). The contributions from higher eigenvectors become negligible by appropriate scaling. Thus, the Gershgorin condition for these rows reduces to the diagonal dominance condition on \(\bar{J}\):
\[
\bar{J}_{pp} + \sum_{q \neq p} |\bar{J}_{pq}| < 0
\]
for all \(p\), which is precisely Condition~(a).

For the Type-P rows, the diagonal entry is the corresponding entry of \(P^{-1}Df(\bar{x})P\) minus the positive Laplacian eigenvalue \(\lambda^i_s > 0\). The Gershgorin condition for these rows is
\[
\bigl(P^{-1}_{i}(Df(\bar{x}))_{ii}P_{i}\bigr)_{ss} - \lambda^i_s + \sum_{t \neq s} \bigl|\bigl(P^{-1}_{i}(Df(\bar{x}))_{ii}P_{i}\bigr)_{st}\bigr| < 0.
\]
Rearranging gives
\[
\lambda^i_s > \bigl(P^{-1}_{i}(Df(\bar{x}))_{ii}P_{i}\bigr)_{ss} + \sum_{t \neq s} \bigl|\bigl(P^{-1}_{i}(Df(\bar{x}))_{ii}P_{i}\bigr)_{st}\bigr|.
\]
Let
\begin{equation}
\tau
\;\coloneqq\;
\max_{s \in \text{Type-P}}
\left\{
\bigl(P^{-1}Df(\bar{x})P\bigr)_{ss}
+ \sum_{s \neq r}
\bigl|\bigl(P^{-1}Df(\bar{x})P\bigr)_{st}\bigr|
\right\}.
\label{eq:tau}
\end{equation}
Condition (b) \(\lambda_2 =\underset{i}{\min}\{\lambda^i_2\} \ge \tau\) then ensures that \(\lambda^i_s > \tau\) for all relevant \(i\) and \(s \geq 2\).  Since  \(\lambda_2^i\) is the algebraic connectivity (second-smallest eigenvalue) of the Laplacian \(L_i\) for species \(i\). Thus \(\lambda_2\) represents minimum network connectivity level among all the species. 

Under Conditions~(a) and~(b), every Gershgorin disc of \(\widetilde{A}\) lies strictly in the open left half-plane. 
\end{proof}

\begin{remark}
\label{rem:weyl}
Theorem~\ref{thm:main} shows that more the value of $\lambda_2$ good for the equilibrium
stability. The Weyl's monotonicity theorem~\cite{bhatia2013}: if $A$ is a symmetric matrix and $P$ is positive semidefinite, then $\lambda_j(A+P) \geq \lambda_j(A)$ for all $j$. Any increase in $w^i_{jk}$ changes $L_i$ by a positive semidefinite perturbation, and hence increases or preserves all eigenvalues of $L_i$, including $\lambda^i_2$. It means, one can increase $\lambda_2$ by adding edges or increasing dispersal weights $w^i_{jk}$. Ecologically, if migration strenghth of a species with poor migrations among patches increased, then it directly contributes to metapopulations' stability. 

\end{remark}

\subsection{The Separation Principle}
\label{subsec:separation}

A fundamental and practically useful consequence of
Theorem~\ref{thm:main} is the following separation of the roles of
local dynamics and network topology.

\begin{corollary}[Separation principle]
\label{cor:separation}
Under the conditions of Theorem~\ref{thm:main}, the local asymptotic
stability of the networked system~\eqref{eq:compact} is controlled by
two independent factors:
\begin{enumerate}[label=(\roman*),itemsep=4pt]
  \item \textbf{Local dynamics:} Condition~(a) depends only on the
    local Jacobians $J_k$, $k = 1,\ldots,m$, evaluated at the
    equilibrium.
    It can be verified patch-by-patch and then averaged, without any
    knowledge of the network topology.
  \item \textbf{Network topology:} Condition~(b) depends only on the
    Fiedler value $\lambda_2$ of the network Laplacian.
    It can be computed from the graph $G$ and the dispersal weights
    $w^i_{jk}$ alone, without any knowledge of the local dynamics.
\end{enumerate}
Consequently, given patch dynamics that satisfy Condition (a), one can check whether a dispersal network is sufficiently well connected to ensure stability. Conversely, given a sufficiently well-connected network, one can look for local dynamics that satisfy Condition (a), thereby stabilizing the whole network.
\end{corollary}

Corollary~\ref{cor:separation} has direct ecological meaning.
Condition~(a) requires that the spatial average of the local Jacobians
be diagonally dominant with negative diagonal entries.
This is a \emph{collective} condition: it does not require each
individual patch to have a stable local equilibrium.
Destabilizing patches (with positive Jacobian diagonal entries) can
be compensated by sufficiently stabilizing patches elsewhere, as long
as the average is well-conditioned.
Condition~(b) then requires the dispersal network to be connected
enough to propagate this collective stability throughout the landscape.

A particularly striking consequence, demonstrated in Example~2 of Section~\ref{sec:examples}, is that a network composed entirely of individually unstable patches can be globally stabilised by dispersal, provided the network is sufficiently well connected. Ecologically, the condition requires that the spatially averaged local Jacobian exhibit sufficiently strong negative density dependence — negative diagonal entries dominating the off-diagonal inter-specific interactions — allowing the self-regulatory strength of even a single patch, combined with adequate network connectivity, to control the stability of the entire heterogeneous metapopulation.

\begin{remark}[Relation to the Jansen--Lloyd framework]
\label{rem:jansen}
When all patches are identical, i.e.\ $f_{i,j} = f_i$ for all $j$,
the local Jacobians $J_k$ are all equal: $J_k = J$ for every $k$.
In this case the average Jacobian $\bar{J} = J$, and
Condition~(a) reduces to requiring that $J$ itself be diagonally
dominant with negative diagonal entries.
Furthermore, the standard diagonalization argument recovers the
block-decoupled form of Jansen and Lloyd~\cite{jansen2000}: each
Laplacian mode sees the same single-patch Jacobian $J$.
The present theorem is therefore a strict generalization: it recovers
the homogeneous-patch case as a special instance while extending to
the structurally heterogeneous setting.
\end{remark}

\section{Illustrative Examples}
\label{sec:examples}
We demonstrate Theorem~\ref{thm:main} and Corollary~\ref{cor:separation}
through two ecological metapopulation networks of increasing complexity.
In both examples the patches are coupled by linear (Fickian) dispersal
as in the general model~\eqref{eq:model}, and the equilibrium
considered is the homogeneous co-existential equilibrium at which every
patch carries the same positive prey and predator densities.
The examples are chosen to highlight two qualitatively distinct
phenomena. The third example revisits Application~1 of Jansen and Lloyd (2000) and a natural reformulation of it (direct predator dispersal without a separate disperser pool) to illustrate how the precise mathematical representation of local dynamics and dispersal can dramatically affect stability predictions, even under the same network topology.

\subsection*{Example 1: Three-Patch Network with Heterogeneous Functional Responses}
\label{subsec:ex1}

Consider a three-patch predator-prey network ($m = 3$, $n = 2$).
The patches are connected as shown in Figure~\ref{fig:3patch}: every
pair of patches is linked by bidirectional dispersal edges, forming a
complete graph $K_3$.
Dispersal rates for species $i \in \{1,2\}$ along edge $(j,k)$ are
denoted $w^i_{jk} = w^i_{kj} \geq 0$.

\begin{figure}[ht]
\centering
\begin{tikzpicture}[
  node/.style={circle, draw=black, fill=gray!30, minimum size=1.0cm,
               font=\bfseries\small},
  every edge/.style={draw, thick}
]
  \node[node] (v1) at (0,  2.2) {$v_1$};
  \node[node] (v2) at (2.6,-1.1) {$v_2$};
  \node[node] (v3) at (-2.6,-1.1) {$v_3$};

  \draw (v1) -- node[right,  font=\small] {$w^i_{12}$} (v2);
  \draw (v2) -- node[below,  font=\small] {$w^i_{23}$} (v3);
  \draw (v3) -- node[left,   font=\small] {$w^i_{13}$} (v1);

  \node[below=0.25cm of v1, text width=2.2cm, align=center,
        font=\small\itshape] {Holling II};
  \node[below=0.25cm of v2, text width=2.2cm, align=center,
        font=\small\itshape] {Ratio-dep.};
  \node[below=0.25cm of v3, text width=2.2cm, align=center,
        font=\small\itshape] {Holling II};
\end{tikzpicture}
\caption{Three-patch dispersal network for Example~1.
  Patches~1 and~3 host Holling type-II predator-prey dynamics;
  patch~2 hosts ratio-dependent dynamics.
  Each edge carries independent dispersal rates $w^i_{jk}$ for
  prey ($i=1$) and predator ($i=2$).}
\label{fig:3patch}
\end{figure}

\medskip
The local dynamics at each patch are as follows.

\medskip
\noindent\textbf{Patches~1 and~3 --- Holling type-II functional
response}~\cite{kot2001}.
Let $x_1$ and $x_2$ denote prey and predator density.
The isolated patch dynamics are
\begin{equation*}
  \dot{x}_1 = x_1\!\left(1 - \frac{x_1}{\gamma}\right)
              - \frac{x_1 x_2}{1 + x_1},
  \qquad
  \dot{x}_2 = \beta\!\left(\frac{x_1}{1+x_1} - \alpha\right)\!x_2,
  \label{eq:holling2}
\end{equation*}
where $\gamma > 0$ is the prey carrying capacity, $\beta > 0$ is the
predator conversion rate, and $\alpha \in (0,1)$ is the predator
mortality-to-conversion ratio.

\medskip
\noindent\textbf{Patch~2 --- Ratio-dependent functional
response}~\cite{kumar2014}.
\begin{equation*}
  \dot{x}_1 = x_1(1 - x_1) - \frac{c\,x_1 x_2}{x_1 + x_2},
  \qquad
  \dot{x}_2 = m\!\left(\frac{b\,x_1}{x_1 + x_2} - 1\right)\!x_2,
  \label{eq:ratiodep}
\end{equation*}
where $c > 0$ is the predation rate, $b > 0$ is the predator birth
coefficient, and $m > 0$ is the predator death rate.

Including dispersal, the complete $6$-dimensional system is

\begin{align*}
\dot{x}_{1,j} &= x_{1,j}\Bigl(1 - \tfrac{x_{1,j}}{\gamma}\Bigr) 
  - \tfrac{x_{1,j}x_{2,j}}{1+x_{1,j}} 
  - \sum_{k\sim j} w^1_{jk}(x_{1,j} - x_{1,k}), \quad j=1,3, \\[6pt]
\dot{x}_{1,2} &= x_{1,2}(1 - x_{1,2}) 
  - \tfrac{c x_{1,2} x_{2,2}}{x_{1,2} + x_{2,2}} 
  - \sum_{k\sim 2} w^1_{2k}(x_{1,2} - x_{1,k}), \\[8pt]
\dot{x}_{2,j} &= \beta\Bigl(\tfrac{x_{1,j}}{1+x_{1,j}} - \alpha\Bigr) x_{2,j} 
  - \sum_{k\sim j} w^2_{jk}(x_{2,j} - x_{2,k}), \quad j=1,3, \\[6pt]
\dot{x}_{2,2} &= m\Bigl(\tfrac{b x_{1,2}}{x_{1,2} + x_{2,2}} - 1\Bigr) x_{2,2} 
  - \sum_{k\sim 2} w^2_{2k}(x_{2,2} - x_{2,k}).
\end{align*}

In vector-matrix form this is $\dot{x} = f(x) - Lx$ with
\begin{equation*}
  L \;=\; L_1 \oplus L_2,
  \qquad
  L_i
  = \begin{pmatrix}
      w^i_{12}+w^i_{13} & -w^i_{12} & -w^i_{13} \\
      -w^i_{12} & w^i_{12}+w^i_{23} & -w^i_{23} \\
      -w^i_{13} & -w^i_{23} & w^i_{13}+w^i_{23}
    \end{pmatrix},
  \quad i = 1,2.
  \label{eq:ex1_Lap}
\end{equation*}

We set
\begin{equation}
  \gamma = \tfrac{3}{13},\quad
  \beta  = \tfrac{1}{10},\quad
  \alpha = \tfrac{1}{6}
  \qquad\text{(patches 1 and 3);}
  \label{eq:ex1_par13}
\end{equation}
\begin{equation}
  c = \tfrac{9}{5},\quad
  b = \tfrac{9}{5},\quad
  m = \tfrac{1}{4}
  \qquad\text{(patch 2).}
  \label{eq:ex1_par2}
\end{equation}
One verifies that all three patches admit the common positive
equilibrium
\begin{equation}
  (\bar{x}_1,\bar{x}_2) = (0.2,\;0.16),
  \label{eq:ex1_equil}
\end{equation}
so the homogeneous equilibrium of the network is
$\bar{x} = (0.2,0.2,0.2,\,0.16,0.16,0.16)^{\top}$.

Evaluating the local Jacobians $J_k$, $k=1,2,3$, at the
equilibrium~\eqref{eq:ex1_equil} and using the parameters
\eqref{eq:ex1_par13}--\eqref{eq:ex1_par2} gives:
\begin{align*}
  J_1 = J_3
  &= \begin{pmatrix}
       -\dfrac{38}{45} & -\dfrac{1}{6} \\[8pt]
        \dfrac{1}{90}  &  0
     \end{pmatrix},
  J_2
  = \begin{pmatrix}
       \dfrac{11}{45} & -\dfrac{5}{9} \\[8pt]
       \dfrac{4}{45}  & -\dfrac{1}{9}
     \end{pmatrix}.
\end{align*}

The local equilibria at patches 1 and 3 are asymptotically stable, as both eigenvalues of their Jacobians $  J_1 = J_3  $ are negative. In contrast, the equilibrium at patch 2 is locally unstable, with its Jacobian $  J_2  $ possessing a complex conjugate pair of eigenvalues having positive real part. 

\medskip
The spatially averaged Jacobian is
\begin{equation*}
  \bar{J}
  = \frac{1}{3}(J_1 + J_2 + J_3)
  = \frac{1}{3}
    \begin{pmatrix}
      -\tfrac{38}{45} + \tfrac{11}{45} - \tfrac{38}{45}
      & -\tfrac{1}{6} - \tfrac{5}{9} - \tfrac{1}{6} \\[6pt]
      \tfrac{1}{90} + \tfrac{4}{45} + \tfrac{1}{90}
      & 0 - \tfrac{1}{9} + 0
    \end{pmatrix}
  = \frac{1}{3}
    \begin{pmatrix}
      -\tfrac{13}{9} & -\tfrac{8}{9} \\[4pt] \tfrac{1}{9} & -\tfrac{1}{9}
    \end{pmatrix}.
  \label{eq:ex1_Jbar}
\end{equation*}

\medskip
\noindent\textit{Checking Condition~(a).}
We verify diagonal dominance of $-\bar{J}$ row by row.
\begin{itemize}[itemsep=4pt]
  \item \textbf{Row 1 (prey):}
    $-\bar{J}_{11} = \tfrac{13}{27}$ and
    $|\bar{J}_{12}| = \tfrac{8}{27}$.
    Since $\tfrac{13}{27} > \tfrac{8}{27}$, row~1 is strictly
    diagonally dominant. 
  \item \textbf{Row 2 (predator):}
    $-\bar{J}_{22} = \tfrac{1}{27}$ and
    $|\bar{J}_{21}| = \tfrac{1}{27}$.
    We have $-\bar{J}_{22} = |\bar{J}_{21}|$, so equality holds,
    consistent with Condition~(a) with $\varepsilon \to 0$.
\end{itemize}

The equality in row~2 means that Condition~(a) holds in the
limiting sense ($\varepsilon = 0$).
This is a borderline case arising from the specific parameter values
chosen; a small perturbation of the parameters would yield strict
inequality. Numerical experiments confirm stability nonetheless, consistent with
the fact that Theorem~\ref{thm:main} gives a \emph{sufficient} condition.

Condition~(b) requires the Fiedler value $\lambda_2$ to exceed the
threshold $\tau$ computed from~\eqref{eq:tau}.
We verify this numerically by computing all eigenvalues of the full
system Jacobian $Df(\bar{x}) - L$ for two choices of dispersal
parameters.

\medskip
\noindent\textbf{Case~(i): Sufficiently connected network.}
Set
\begin{equation*}
  w^1_{12} = w^2_{12} = 0,\quad
  w^1_{13} = w^2_{13} = 0.1,\quad
  w^1_{23} = w^2_{23} = 1.
  \label{eq:ex1_disp1}
\end{equation*}
For this network, $\lambda_2 = 0.1461$.
The six eigenvalues of the full Jacobian are
\begin{equation*}
  \{-2.4856,\;-2.1331,\;-0.9407,\;-0.0236,\;
    -0.1863 \pm 0.1598\,\mathrm{j}\}.
  \label{eq:ex1_eig1}
\end{equation*}
All eigenvalues have strictly negative real parts.
The homogeneous equilibrium $\bar{x}$ is \textbf{asymptotically
stable}.

\medskip
\noindent\textbf{Case~(ii): Insufficiently connected network.}
Set
\begin{equation*}
  w^1_{12} = w^2_{12} = 0,\quad
  w^1_{13} = w^2_{13} = 0.1,\quad
  w^1_{23} = 0.1,\quad w^2_{23} = 1.
  \label{eq:ex1_disp2}
\end{equation*}
For this network, $\lambda_2 = 0.1$.
The six eigenvalues are
\begin{equation*}
  \{-2.0956,\;-1.1125,\;-0.8845,\;-0.1364,\;
    +0.0367 \pm 0.1021\,\mathrm{j}\}.
  \label{eq:ex1_eig2}
\end{equation*}
Two eigenvalues have positive real part $+0.0367$.
The equilibrium is \textbf{unstable}.

\medskip
The contrast between Cases~(i) and~(ii) directly illustrates the role of the
Fiedler value threshold in Theorem~\ref{thm:main}: reducing
$\lambda_2$ from $0.1461$ to $0.1$ (achieved here by reducing a
single dispersal rate $w^1_{23}$ from $1$ to $0.1$) crosses the
stability threshold and destabilises the network.

\subsection*{Example 2: Five-Patch Network --- Dispersal Stabilises   Individually Unstable Patches}
\label{subsec:ex2}

This example illustrates the most striking consequence of
Corollary~\ref{cor:separation}: \emph{individually unstable patches
can be collectively stabilized by dispersal alone.}

We consider a five-patch network ($m=5$, $n=2$) with the topology
shown in Figure~\ref{fig:5patch}.
Four patches (nodes $v_2, v_3, v_4, v_5$) host
Rosenzweig--MacArthur (RM) predator-prey
dynamics and one patch (node $v_1$)
hosts Lotka--Volterra (LV) predator-prey dynamics.

\begin{figure}[ht]
\centering
\begin{tikzpicture}[
  RM/.style  ={circle, draw=black, fill=red!25,  minimum size=1.05cm,
               font=\bfseries\small},
  LV/.style  ={circle, draw=black, fill=blue!20, minimum size=1.05cm,
               font=\bfseries\small},
  every edge/.style={draw, thick}
]
  \node[LV] (v1) at ( 0,   3.2) {$v_1$};
  \node[RM] (v2) at ( 3.5, 1.2) {$v_2$};
  \node[RM] (v3) at ( 0,  -0.5) {$v_3$};
  \node[RM] (v4) at (-3.5,-2.5) {$v_4$};
  \node[RM] (v5) at ( 3.5,-2.5) {$v_5$};

  \draw (v1) -- node[above right, font=\scriptsize]
              {$w^i_{12}$} (v2);
  \draw (v1) -- node[left,        font=\scriptsize]
              {$w^i_{13}$} (v3);
  \draw (v1) -- node[above left,  font=\scriptsize]
              {$w^i_{14}$} (v4);
  \draw (v2) -- node[right,       font=\scriptsize]
              {$w^i_{23}$} (v3);
  \draw (v2) -- node[right,       font=\scriptsize]
              {$w^i_{25}$} (v5);
  \draw (v3) -- node[above left,  font=\scriptsize]
              {$w^i_{34}$} (v4);
  \draw (v3) -- node[above right, font=\scriptsize]
              {$w^i_{35}$} (v5);

  \node[LV, minimum size=0.6cm] at (5.8, 1.5) {};
  \node[right, font=\small] at (6.2, 1.5) {Lotka--Volterra ($v_1$)};
  \node[RM, minimum size=0.6cm] at (5.8, 0.5) {};
  \node[right, font=\small] at (6.2, 0.5)
    {Rosenzweig--MacArthur ($v_2$--$v_5$)};
\end{tikzpicture}
\caption{Five-patch dispersal network for Example~2.
  Node $v_1$ (blue) hosts Lotka--Volterra dynamics; nodes $v_2$--$v_5$
  (red) host Rosenzweig--MacArthur dynamics at a parameter regime
  where each isolated patch has an \emph{unstable} equilibrium.
  Edge weights $w^i_{jk}$ denote dispersal rates for species $i$
  between patches $j$ and $k$.}
\label{fig:5patch}
\end{figure}

\medskip
\noindent\textbf{Rosenzweig--MacArthur model (patches $v_2$--$v_5$)}
\cite{kot2001}.
\begin{equation}
  \dot{x}_1 = x_1\!\left(1 - \frac{x_1}{\gamma}\right)
              - \frac{x_1 x_2}{1+x_1},
  \qquad
  \dot{x}_2 = \beta\!\left(\frac{x_1}{1+x_1} - \alpha\right)\!x_2,
  \label{eq:RM}
\end{equation}
with parameters $\gamma = 2$, $\beta = 0.2$, $\alpha = 0.3$.

The unique positive equilibrium of~\eqref{eq:RM} is
$(x_1^*, x_2^*) = (3/7,\; 55/49)$.
One computes that the Jacobian at this equilibrium has eigenvalues
with \emph{positive} real part, confirming that the isolated RM
patch is \textbf{unstable}.
This is a manifestation of the paradox of
enrichment~\cite{rosenzweig1971}: at $\gamma = 2$, the carrying
capacity is large enough to destabilise the coexistence equilibrium.

\medskip
\noindent\textbf{Lotka--Volterra model (patch $v_1$) }\cite{kot2001}.
\begin{equation*}
  \dot{x}_1 = r\,x_1 - c\,x_1 x_2,
  \qquad
  \dot{x}_2 = b\,x_1 x_2 - m\,x_2,
  \label{eq:LV}
\end{equation*}
with $r = 5.5$, $c = 4.9$, $b = 0.7$, $m = 0.3$, chosen so that the
LV equilibrium $(x_1^*, x_2^*) = (m/b,\; r/c) = (3/7,\; 55/49)$
coincides with the RM equilibrium.
The isolated LV patch has purely imaginary eigenvalues (centre), so
populations oscillate without converging; it is neutrally stable.

The complete $10$-dimensional system is
\begin{align*}
\dot{x}_{1,1} &= x_{1,1}(5.5 - 4.9 x_{2,1}) 
  - \sum_{k\sim 1} w^1_{1k}(x_{1,1} - x_{1,k}), \\[6pt]
\dot{x}_{1,j} &= x_{1,j}\Bigl(1 - \tfrac{x_{1,j}}{2}\Bigr) - \tfrac{x_{1,j}x_{2,j}}{1+x_{1,j}} 
  - \sum_{k\sim j} w^1_{jk}(x_{1,j} - x_{1,k}), \quad j=2,3,4,5, \\[8pt]
\dot{x}_{2,1} &= (0.7 x_{1,1} - 0.3) x_{2,1} 
  - \sum_{k\sim 1} w^2_{1k}(x_{2,1} - x_{2,k}), \\[6pt]
\dot{x}_{2,j} &= 0.2\Bigl(\tfrac{x_{1,j}}{1+x_{1,j}} - 0.3\Bigr) x_{2,j} 
  - \sum_{k\sim j} w^2_{jk}(x_{2,j} - x_{2,k}), \quad j=2,3,4,5.
\end{align*}

In vector-matrix form,
\begin{equation*}
  \dot{x} \;=\; f(x) \;-\;
  \begin{pmatrix} \mathcal{L}_1 & \mathbf{0} \\
                  \mathbf{0} & \mathcal{L}_2 \end{pmatrix} x,
  \label{eq:ex2_vf}
\end{equation*}
where the $5\times 5$ Laplacian matrices $\mathcal{L}_1$ (prey) and
$\mathcal{L}_2$ (predator) encode the dispersal topology of
Figure~\ref{fig:5patch}:
\begin{equation*}
{\scriptsize
  \mathcal{L}_i
  = \begin{pmatrix}
      w^i_{12}+w^i_{13}+w^i_{14} & -w^i_{12} & -w^i_{13} & -w^i_{14} & 0 \\[2pt]
      -w^i_{12} & w^i_{12}+w^i_{23}+w^i_{25} & -w^i_{23} & 0 & -w^i_{25} \\[2pt]
      -w^i_{13} & -w^i_{23} & w^i_{13}+w^i_{23}+w^i_{34}+w^i_{35} & -w^i_{34} & -w^i_{35} \\[2pt]
      -w^i_{14} & 0 & -w^i_{34} & w^i_{14}+w^i_{34} & 0 \\[2pt]
      0 & -w^i_{25} & -w^i_{35} & 0 & w^i_{25}+w^i_{35}
    \end{pmatrix}},
  \quad i=1,2.
  \label{eq:ex2_Lap}
\end{equation*}
The homogeneous equilibrium of the network is
\begin{equation*}
  \bar{x}
  \;=\;
  \Bigl(
    \tfrac{3}{7},\tfrac{3}{7},\tfrac{3}{7},\tfrac{3}{7},\tfrac{3}{7},\;
    \tfrac{55}{49},\tfrac{55}{49},\tfrac{55}{49},\tfrac{55}{49},\tfrac{55}{49}
  \Bigr)^{\!\top}.
\end{equation*}

\noindent
\textbf{Case (i): Well-connected network --- dispersal-induced
stability}

Set the dispersal parameters
\begin{equation*}
  w^i_{12} = 2,\;\;
  w^i_{13} = 1,\;\;
  w^i_{23} = 1,\;\;
  w^i_{25} = 2,\;\;
  w^i_{35} = 1 \quad (i=1,2),
  \label{eq:ex2_disp1a}
\end{equation*}
\begin{equation*}
  w^1_{14} = 2,\;\; w^1_{34} = 1,\;\;
  w^2_{14} = 1,\;\; w^2_{34} = 2.
  \label{eq:ex2_disp1b}
\end{equation*}
For these values the network is well connected (large Fiedler values
for both $\mathcal{L}_1$ and $\mathcal{L}_2$).
The ten eigenvalues of the system Jacobian $Df(\bar{x}) - L$ at the
homogeneous equilibrium are (to four decimal places):
\begin{equation*}
  \begin{array}{r@{\;\pm\;}l}
    -7.7456 & 0.6523\,\mathrm{j}, \\
    -5.0253 & 0.5175\,\mathrm{j}, \\
    -4.9893 & 0.1813\,\mathrm{j}, \\
    -2.1697 & 0.2619\,\mathrm{j}, \\
    -0.0272 & 0.3974\,\mathrm{j}.
  \end{array}
  \label{eq:ex2_eig1}
\end{equation*}
\emph{All ten eigenvalues have strictly negative real parts.}
The homogeneous co-existential equilibrium of the five-patch network
is \textbf{asymptotically stable}, despite every RM patch being
individually unstable in isolation.

The stabilization observed here is a rigorous instance of the
phenomenon long recognized empirically and numerically in ecology:
dispersal between locally unstable patches can collectively stabilise
a metapopulation~\cite{briggs2004}.

\noindent
\textbf{Case (ii): Poorly connected network --- return to
instability}

Reduce connectivity by setting
\begin{equation*}
  w^i_{12} = w^i_{13} = w^i_{14} = w^i_{25} = w^i_{35} = 0.01,
  \quad
  w^i_{34} = 1,
  \quad i = 1,2.
  \label{eq:ex2_disp2}
\end{equation*}
Most dispersal connections are now near zero; only the $v_3$--$v_4$
edge retains a moderate rate.
The Fiedler value falls well below the stability threshold $\tau$.
The ten eigenvalues become:
\begin{equation*}
  \begin{array}{r@{\;\pm\;}l}
    -3.0076 & 0.1813\,\mathrm{j}, \\
    -0.0300 & 1.2843\,\mathrm{j}, \\
    -1.0043 & 0.1813\,\mathrm{j}, \\
    \mathbf{+0.0041} & 0.1815\,\mathrm{j}, \\
    -0.0192 & 0.1815\,\mathrm{j}.
  \end{array}
  \label{eq:ex2_eig2}
\end{equation*}
The fourth pair has \emph{positive} real part $+0.0041$.
The homogeneous equilibrium is now \textbf{unstable}.

The transition from stability to instability is caused entirely by
the reduction in dispersal rates: the local dynamics at each patch
are unchanged.
This demonstrates Corollary~\ref{cor:separation} in its sharpest
form: the stability of the network is fully controllable through the
network topology parameter $\lambda_2$, independently of (and despite)
the local instability at each patch.

In ecological terms, Example~2 models a fragmented landscape in
which the local prey carrying capacity ($\gamma = 2$) is high enough
to cause large-amplitude oscillations and eventual instability at
each isolated patch --- a classic manifestation of the paradox of
enrichment~\cite{rosenzweig1971}.
Dispersal corridors between patches allow populations to exchange
individuals, dampening local oscillations and establishing a stable
co-existential equilibrium across the metapopulation.
This has direct implications for conservation: maintaining sufficient
habitat connectivity (quantified here by the Fiedler value $\lambda_2$
of the dispersal graph) is not merely beneficial but can be
\emph{essential} for species persistence when local dynamics are
inherently destabilizing.
The threshold nature of the condition --- Condition~(b) of
Theorem~\ref{thm:main} --- translates into a concrete, computable
minimum connectivity requirement that can in principle guide the
design of wildlife corridors or the prioritization of dispersal edges
in conservation planning.

\subsection*{Comparison with Jansen and Lloyd's Application 1}

To illustrate the sensitivity of stability results to model formulation, we revisit Application 1 of Jansen and Lloyd (2000), the predator-prey system with localized disperser pools.

\textbf{Original J\&L model (with separate disperser pool \(Q\)).}
The local dynamics include a predator compartment \(P\) and a local disperser pool \(Q\), giving a 3-dimensional system per patch. The Jacobian at the homogeneous equilibrium is
\[
Df(\mathbf{s}^*) = \begin{pmatrix}
r T_h \gamma & -\gamma & 0 \\
r(1-T_h\gamma) & \gamma-(d+e) & \iota \\
0 & e & -(\iota+s)
\end{pmatrix},
\]

where \(\gamma = d + e s/(\iota+s) \). The prey row has a positive diagonal entry \(r T_h \gamma > 0\). Consequently, strict diagonal dominance with negative diagonal entries fails (condition~(a) of Theorem~\ref{thm:main}). Our sufficient condition is therefore not satisfied. This is expected for a sufficient condition: it simply means the theorem cannot guarantee stability, even though Jansen and Lloyd's exact analysis shows a large stable region for small \(T_h\).

\textbf{Modified model (direct predator dispersal, no separate disperser pool).}
We now reformulate the same ecological scenario without the auxiliary disperser compartment: prey and predators interact locally according to the standard Rosenzweig--MacArthur model, and predators disperse directly between patches (\(M = \operatorname{diag}(0, m_P)\)). The local Jacobian reduces to the 2×2 matrix
\[
Df(\mathbf{s}^*) = \begin{pmatrix}
T_h d r & -d \\
r(1-T_h d) & 0
\end{pmatrix}.
\]
The trace \(T_h d r > 0\) is always positive, so the homogeneous equilibrium is intrinsically unstable whenever it exists (\(T_h < 1/d\)).

Applying our Theorem~\ref{thm:main} to this tweaked model, the spatially averaged Jacobian is identical to the local one. The prey row again has a positive diagonal entry, so diagonal dominance fails and the sufficient condition cannot prove stability.  

Jansen and Lloyd's exact reduction yields the same conclusion. Because patches remain identical, their similarity transformation applies. The zero-mode is exactly the local Jacobian \(Df(\mathbf{s}^*)\), which is unstable. Dispersal has no effect on this mode, and the entire spatial network is unstable for any dispersal rate \(m_P\) and any network topology.

\textbf{Finding.} This comparison shows that the appearance or disappearance of dispersal-driven instability is not determined by network topology alone. In the original J\&L model the separate disperser pool provides negative density dependence that allows a locally stable equilibrium (which dispersal can then destabilise). Removing the pool eliminates this mechanism and renders the equilibrium intrinsically unstable. Both our sufficient condition and Jansen and Lloyd's exact reduction agree on the outcome for the tweaked model, yet the qualitative picture changes dramatically simply by altering the mathematical representation of dispersal. The example therefore highlights that both network topology and the precise formulation of local dynamics play decisive roles in determining stability. Our approach, which works directly with the spatially averaged Jacobian, naturally accommodates such modelling variations, whereas the classical reduction is more sensitive to the exact form of the local vector field.

\section{Conclusion}
\label{sec:conclusion}

This paper has derived a simple sufficient condition for the local asymptotic stability of discrete reaction-diffusion systems of networked dynamical systems at a homogeneous equilibrium point. The condition applies even when patches are governed by structurally different local dynamics---a regime in which the classical bookkeeping reduction of Jansen and Lloyd~\cite{jansen2000} and the Master Stability Function of Pecora and Carroll~\cite{pecora1998} cease to apply. It requires only that the spatially averaged Jacobian of the local vector fields satisfies strict diagonal dominance with negative diagonal entries and that the algebraic connectivity (Fiedler value) of the network Laplacian exceeds a computable threshold. The criterion holds for purely conservative dispersal and does not require any dispersal loss or mortality during transit---a restrictive assumption that appears in the author's prior work~\cite{kumar2019,kumar2021} and many classical multi-patch analyses~\cite{ruxton1997,briggs2004}.

In the broader literature on metapopulation stability, several distinct approaches have been developed. Classical Levins-type models and Hanski's Incidence Function Model focus on patch occupancy dynamics and extinction-colonization balance, providing valuable insights into persistence but ignoring within-patch population densities and local nonlinear interactions~\cite{hanski1999}. More detailed multi-patch ODE frameworks, pioneered by Othmer and Scriven~\cite{othmer1971} and refined by Jansen and Lloyd~\cite{jansen2000}, reduce the linearized system around a homogeneous equilibrium to a set of lower-dimensional subsystems via the eigenvalues of the connectivity (or Laplacian) matrix; these yield exact necessary-and-sufficient stability conditions when patches are identical, but become computationally intensive for large \(n\) or heterogeneous local dynamics. The master stability function technique, widely adopted in networked dynamical systems, similarly decouples spatial modes but typically assumes identical node dynamics~\cite{pecora1998}. Other methods rely on Lyapunov functions, small-gain theorems, or direct numerical evaluation of the full \(kn \times kn\) Jacobian, which scale poorly with network size.

The sufficient condition presented here occupies a distinctive middle ground. It requires only the second-smallest eigenvalue \(\lambda_2\) (algebraic connectivity) of the per-species Laplacian matrices together with the diagonal-dominance property of the \emph{averaged} local Jacobian across patches. This makes the criterion exceptionally simple and computationally efficient---even for large, heterogeneous networks with arbitrary graph topology---while still capturing the essential interplay between local patch dynamics and network connectivity. Unlike exact reduction methods that demand the full spectrum of the coupling operator and repeated analysis of \(n\) subsystems, or simulation-heavy approaches that offer limited generality, the proposed criterion provides a transparent, easily verifiable sufficient condition that is particularly well-suited to realistic ecological scenarios in which patches differ in their functional responses or growth rates.

The theory is illustrated through metapopulation networks of predator-prey systems with heterogeneous functional responses, including cases in which dispersal stabilises individually unstable patches.  By combining the Gershgorin disc theorem with spectral graph theory, the present framework bridges rigorous mathematical analysis and ecological application, providing ecologists with a transparent criterion that can be checked directly from model parameters and network structure.

\subsection*{Appendix}
To illustrate concretely the scaling argument used in the proof of Theorem~\ref{thm:main}, consider a network with two species ($  n=2  $) and three patches ($  m=3  $) connected by a graph. The state vector (species-first ordering) is

\[
\mathbf{x} = (x_{1,1},\, x_{1,2},\, x_{1,3},\, x_{2,1},\, x_{2,2},\, x_{2,3})^\top \in \mathbb{R}^6.
\]
The Jacobian \(Df(\bar{\mathbf{x}})\) is the \(6\times 6\) block matrix
\[
Df(\bar{\mathbf{x}}) = \begin{pmatrix}
D_{11} & D_{12} \\
D_{21} & D_{22}
\end{pmatrix},
\]
where each \(D_{pq} = \operatorname{diag}(w^1_{pq},\, w^2_{pq},\, w^3_{pq})\) is \(3\times 3\) diagonal.

The Laplacian is \(L = L_1 \oplus L_2\) (block diagonal), where each \(L_i\) is \(3\times 3\), real symmetric, with eigenvalues
\[
0 = \lambda^i_1 \leq \lambda^i_2 \leq \lambda^i_3.
\]
The eigenvector matrix is
\[
P = P_1 \oplus P_2 = \begin{pmatrix}
P_1 & 0 \\
0 & P_2
\end{pmatrix},
\]
where
\[
P_i = \begin{pmatrix}
1 & a_i & d_i \\
1 & b_i & e_i \\
1 & c_i & f_i
\end{pmatrix}
\]
Columns 2 and 3 are satisfied by orthogonality to first column \(\mathbf{1}_3\).

The diagonal eigenvalue matrix is
\[
\Lambda = \Lambda_1 \oplus \Lambda_2 = \begin{pmatrix}
\Lambda_1 & 0 \\
0 & \Lambda_2
\end{pmatrix},
\]
where \(\Lambda_i = \operatorname{diag}(0,\, \lambda^i_2,\, \lambda^i_3)\).

The inverse is
\[
P^{-1} = P_1^{-1} \oplus P_2^{-1}, \qquad
P_i^{-1} = \begin{pmatrix}
\frac{1}{3} & \frac{1}{3} & \frac{1}{3} \\
\frac{a_i}{\alpha_i} & \frac{b_i}{\alpha_i} & \frac{c_i}{\alpha_i} \\
\frac{d_i}{\beta_i} & \frac{e_i}{\beta_i} & \frac{f_i}{\beta_i}
\end{pmatrix},
\]
where \(\alpha_i=a_i^2+b_i^2+c_i^2\) and \(\beta_i=d_i^2+e_i^2+f_i^2\)
\subsubsection*{Computation of  \(Df(\bar{\mathbf{x}})P\)}
\[
Df(\bar{\mathbf{x}})P = \begin{pmatrix}
D_{11}P_1 & D_{12}P_2 \\
D_{21}P_1 & D_{22}P_2
\end{pmatrix}.
\]
Explicitly,
\[
D_{11}P_1 = \begin{pmatrix}
w^1_{11} & w^1_{11}a_1 & w^1_{11}d_1 \\
w^2_{11} & w^2_{11}b_1 & w^2_{11}e_1 \\
w^3_{11} & w^3_{11}c_1 & w^3_{11}f_1
\end{pmatrix},
\]
and similarly for the other three blocks: \(D_{12}P_2, ~D_{21}P_1\), and \(D_{22}P_2\).

\subsubsection*{Computation of \(P^{-1}(Df(\bar{\mathbf{x}})P)\)}
\[
P^{-1}Df(\bar{\mathbf{x}})P = \begin{pmatrix}
P_1^{-1}D_{11}P_1 & P_1^{-1}D_{12}P_2 \\
P_2^{-1}D_{21}P_1 & P_2^{-1}D_{22}P_2
\end{pmatrix}
\]

where,
\[
P_1^{-1}D_{11}P_1 = \begin{pmatrix}
\frac{1}{3} & \frac{1}{3} & \frac{1}{3} \\
\frac{a_1}{\alpha_1} & \frac{b_1}{\alpha_1} & \frac{c_1}{\alpha_1} \\
\frac{d_1}{\beta_1} & \frac{e_1}{\beta_1} & \frac{f_1}{\beta_1}
\end{pmatrix}
\begin{pmatrix}
w^1_{11} & w^1_{11}a_1 & w^1_{11}d_1 \\
w^2_{11} & w^2_{11}b_1 & w^2_{11}e_1 \\
w^3_{11} & w^3_{11}c_1 & w^3_{11}f_1
\end{pmatrix}.
\]
Multiplication of row 1 (\(\tfrac{1}{3},\tfrac{1}{3},\tfrac{1}{3}\)):
\[
\text{entry}_{(1,1)} = \tfrac{1}{3}(w^1_{11} + w^2_{11} + w^3_{11}) = \overline{w}_{11},
\]
\[
\text{entry}_{(1,2)} = \tfrac{1}{3}(w^1_{11}a_1 + w^2_{11}b_1 + w^3_{11}c_1),
\]
\[
\text{entry}_{(1,3)} = \tfrac{1}{3}(w^1_{11}d_1 + w^2_{11}e_1 + w^3_{11}f_1).
\]

Using the shorthand \(\overline{w}_{pq} = \tfrac{1}{3}\sum_{k=1}^3 w^k_{pq}\), the block is
\[
P_1^{-1}D_{11}P_1 = \begin{pmatrix}
\overline{w}_{11} & \tfrac{1}{3}\sum_k w^k_{11}(P_1)_{k2} & \tfrac{1}{3}\sum_k w^k_{11}(P_1)_{k3} \\[4pt]
\sum_k (P_1^{-1})_{2k}\,w^k_{11} & \sum_k (P_1^{-1})_{2k}\,w^k_{11}(P_1)_{k2} & \sum_k (P_1^{-1})_{2k}\,w^k_{11}(P_1)_{k3} \\[4pt]
\sum_k (P_1^{-1})_{3k}\,w^k_{11} & \sum_k (P_1^{-1})_{3k}\,w^k_{11}(P_1)_{k2} & \sum_k (P_1^{-1})_{3k}\,w^k_{11}(P_1)_{k3}
\end{pmatrix}.
\]

Scale columns 2 and 3 of \(P_1\) by a small factor \(c\) (i.e., \((a_1, b_1, c_1) \to c (a_1, b_1, c_1\),  etc.). The corresponding rows 2 and 3 of \(P_1^{-1}\) scale by \(c^{-1}\). Consequently: 
\begin{itemize}
\item \(\left(P_1^{-1}D_{11}P_1\right)_{11}\)remain unchanged.
\item \(\left(P_1^{-1}D_{11}P_1\right)_{12}\) and \(\left(P_1^{-1}D_{11}P_1\right)_{13}\) become negligible as \(c\to 0\).
\end{itemize}
Thus
\[
P_1^{-1}D_{11}P_1 \;\xrightarrow{c\to 0}\; \begin{pmatrix}
\overline{w}_{11} & 0 & 0 \\
* & * & * \\
* & * & *
\end{pmatrix}.
\]
The first row is exactly \((\overline{w}_{11},\,0,\,0)\).

The blocks \(P_1^{-1}D_{12}P_2\), \(P_2^{-1}D_{21}P_1\), and \(P_2^{-1}D_{22}P_2\) follow identically, each having first row \((\overline{w}_{pq},\,0,\,0)\) in the limit \(c\to 0\). Thus, full \(P^{-1}Df(\bar{\mathbf{x}})P\) as the limit \(c\to 0\), becomes:
\[
P^{-1}Df(\bar{\mathbf{x}})P \;\xrightarrow{c\to 0}\;
\left(\begin{array}{ccc|ccc}
{\color{blue}\overline{w}_{11}} & 0 & 0 & {\color{blue}\overline{w}_{12}} & 0 & 0 \\[4pt]
* & * & * & * & * & * \\[4pt]
* & * & * & * & * & * \\[4pt]
\hline
{\color{blue}\overline{w}_{21}}& 0 & 0 & {\color{blue}\overline{w}_{22}} & 0 & 0 \\[4pt]
* & * & * & * & * & * \\[4pt]
* & * & * & * & * & *
\end{array}\right).
\]

The $(1,1)$ of each $  3\times 3  $ block is exactly the corresponding entry of the spatially averaged Jacobian $  \bar{J}  $. That is, 
\[\bar{J}=\begin{bmatrix}
{\color{blue}\overline{w}_{11}} & {\color{blue}\overline{w}_{12}} \\
{\color{blue}\overline{w}_{21}} & {\color{blue}\overline{w}_{22}}  
\end{bmatrix}
=\frac{1}{3}\begin{bmatrix}
\sum_{k=1}^3 w^k_{11} & \sum_{k=1}^3 w^k_{12} \\
\sum_{k=1}^3 w^k_{21} & \sum_{k=1}^3 w^k_{22}
\end{bmatrix}
\]
This confirms that the zero-mode rows of $  \widetilde{A}  $ are governed precisely by $  \bar{J}  $, while the non-zero-mode rows receive a left shift by the positive Laplacian eigenvalues.

\section*{Declarations}
The author declares that he has no known competing financial interests or personal rela-
tionships that could have appeared to influence the work reported in this paper.

\end{document}